\documentstyle[11pt,amssymb]{article}

\def\demo{\noindent{\bf Proof. }}
\def\sqr#1#2{{\vcenter{\hrule height.#2pt
	\hbox{\vrule width.#2pt height#1pt \kern#1pt
		\vrule width.#2pt}
	\hrule height.#2pt}}}
\def\square{\mathchoice\sqr64\sqr64\sqr{4}3\sqr{3}3}
\def\QED{\hfill$\square$}
\def\tratto{\mbox{\rule{2mm}{.2mm}$\;\!$}}

\newtheorem{Theorem}{Theorem}[section]

\newtheorem{Corollary}[Theorem]{Corollary}
\newtheorem{Proposition}[Theorem]{Proposition}
\newtheorem{Remark}[Theorem]{Remark}
\newtheorem{Example}[Theorem]{Example}

\begin{document}
\title{
\hfill\mbox{{\scriptsize \begin{tabular}{cl}
To appear:&J. Algebra
\end{tabular}}}\\
\   \\
{\Huge \bf Links of prime ideals and \\
their Rees algebras }\\
\footnotetext{AMS 1991
{\em Mathematics Subject Classification}. Primary 13H10;
Secondary 13C40,  13D40, 13D45, 13H15.}
}

\author{Alberto Corso \quad\quad
Claudia Polini\thanks{The author was partially supported by NSF grant
STC-91-19999.} \\
\begin{tabular}{c}
\       \\
\       \\
Department of Mathematics, Rutgers University \\
New Brunswick, New Jersey 08903 (USA) \\
{\tt E-mail: corso, polini@math.rutgers.edu} 
\end{tabular}
}
\date{}

\maketitle

\vspace{-0.4in}

\begin{abstract}
\noindent In a previous paper we exhibited the 
somewhat surprising property that most direct links of prime ideals in
Gorenstein rings are equimultiple ideals with reduction number $1$.  
This led to the construction of large families of
Cohen--Macaulay Rees algebras. 
The first goal of this paper is to extend this result to arbitrary 
Cohen--Macaulay rings. The means of the proof are changed since
one cannot depend so heavily on linkage theory. We then 
study the structure of the  Rees algebra of these links, more specifically
we describe their canonical module in sufficient detail to be able to
characterize self--linked prime ideals. In the last section
multiplicity estimates for classes of such ideals are established.   
\end{abstract}

\vfill

%
%
%
%
%
%
%
\newpage

\section{Introduction}
Let $R$ be a Noetherian local ring and let $I$ be one of its ideals.
Two of the most important graded algebras associated to $I$ that
have been extensively studied in the past few years are the {\em Rees
algebra} of $I$
$$
R[It]=\bigoplus_{i\geq 0} I^it^i,
$$
and the {\em associated graded ring} 
$$
{\rm gr}_I(R)=\bigoplus_{i\geq 0} I^i/I^{i+1}.
$$
A recurrent theme is to find out conditions that make them
Cohen--Macaulay rings. A tool that has proved its usefulness 
is the notion of the reduction of an ideal (see \cite{NR}). A {\em
reduction} of $I$ is an ideal $J\subset I$ for which there exists an
integer $n$ such that $I^{n+1}=JI^n$; the least such integer is called 
the {\em reduction number of $I$ with respect to $J$}, in symbols
$r_J(I)$. The reduction number of $I$ is the minimum of $r_J(I)$ 
taken over all the minimal reductions of $I$. 
Phrased in other words, $J$ is a reduction of $I$ if and only if the
morphism $ R[Jt] \hookrightarrow R[It]$ is finite. It seems reasonable
to expect to recover some of the properties of $R[It]$ from those
of $R[Jt]$. It seems even more reasonable to expect better results 
on $R[It]$ when $J$ has an amenable structure, such as being a
complete intersection. 

In the case in which the residue field of the ring is infinite,
the minimal number of generators of a minimal reduction of
$I$ does not depend on the reduction; this number, usually
denoted with $\ell(I)$,  is called the {\em analytic spread} 
of $I$ and is equal to the dimension of the fibre ring of $I$, i.e.
$R[It]\otimes R/{\frak m}$. The relation ${\rm height}(I)\leq \ell(I)
\leq {\rm dim}(R)$ is a natural source of several deviation measures. 
One of them was introduced by S. Huckaba and C. Huneke in
\cite{HuHu}, namely the {\em analytic deviation} of an ideal
$$
{\rm ad}(I)=\ell(I)-{\rm height}(I).
$$
${\rm ad}(I)$ measures how far a minimal reduction of an ideal is from 
being a complete intersection. 

The ideals having analytic deviation zero are called {\em
equimultiple ideals}. New families of such ideals were introduced in
\cite{CPV} for Gorenstein rings by the process of linkage. For more
details on linkage we refer the reader to \cite{PS} while, for the time
being, it will be enough to know that
two proper ideals $I$ and $J$ of height $g$ in a Cohen--Macaulay 
ring $R$ are said to be ({\em directly}) {\em linked}, in symbols 
$I\sim J$, if there exists a regular sequence ${\bf z}=(z_1, ..., z_g)
\subset I\cap J$ such that $J={\bf z}\colon I$ and $I={\bf z}\colon J$. 

To be specific, if ${\frak p}$ is a prime ideal of height $g$ and ${\bf
z}=(z_1, \ldots , z_g)\subset {\frak p}$ is a regular sequence,
\cite{CPV} focused on the direct link $I={\bf z}\colon {\frak p}$, and the
following two conditions played significant roles there:

\begin{itemize}
\item[{\rm (L$_1$)}] $R_{{\frak p}}$ is not a regular local ring;
\item[{\rm (L$_2$)}] $R_{{\frak p}}$ is a regular local ring of
dimension at least $2$ and two elements in the sequence ${\bf z}$ lie
in the symbolic square ${\frak p}^{(2)}$. 
\end{itemize}

\noindent Using these notations \cite[Theorem 2.1, Theorem 3.1,
Corollary 3.2]{CPV}:

\medskip

\noindent{\bf Theorem A\ } {\em Let $R$ be a Cohen--Macaulay ring,
${\frak p}$  a prime ideal of codimension $g$, and let ${\bf z}=(z_1,
\ldots, z_g)\subset {\frak p}$ be a regular sequence. Set $J=({\bf z})$
and $I=J\colon {\frak p}$. Suppose that $R_{\frak p}$ is a Gorenstein
ring. Then $I$ is an equimultiple ideal with reduction number one,
more precisely, 
$$
I^2=J I,
$$
if either condition {\rm L}${}_1$ or {\rm L}${}_2$ holds.}

\medskip

\noindent{\bf Theorem B\ } {\em
Let $R$ be a Gorenstein local ring, ${\frak p}$ a Cohen--Macaulay
prime ideal of codimension $g$, $J= (z_1,...,z_g)$ a complete
intersection, and set $I=J\colon{\frak p}$. If either condition {\rm
L}${}_1$ or {\rm L}${}_2$ holds then the associated graded ring ${\rm
gr}_I(R)$ is Cohen--Macaulay.} 

\medskip

\noindent{\bf Corollary C\ } {\em If in {\rm Theorem B} we have ${\rm
dim}(R_{\frak p}) \geq 2$ {\rm (}which is automatic in the second
case{\rm )}, then the Rees algebra of $I$ is Cohen--Macaulay. }

\medskip

As a result, the process described above gives us a powerful
tool for the production of Cohen--Macaulay Rees algebras.

\medskip

It is now straightforward to describe our first main result, the
extension of Theorem A from Gorenstein rings to general
Cohen--Macaulay rings. 

\medskip

\noindent {\bf Theorem \ref{main}\ }{\em
Let $(R,{\frak m})$ be a Cohen--Macaulay local ring and let ${\frak
p}$ be a prime ideal of $R$ of height $g$ such that $R_{\frak p}$ is
not a regular local ring. If $J=(z_1, ..., z_g)$ denotes a regular
sequence inside ${\frak p}$, then the link $I=J\colon {\frak p}$ is
equimultiple with reduction number $1$. }

\medskip

This will be proved in Section $2$.
The approach is fully ideal--theoretic, in contrast to \cite{CPV},
when we could appeal to basic properties of linkage theory. The fact
that we are dealing with direct links of prime ideals still permits certain
calculations to go through.

To benefit however from this close relationship between the link and
its reduction to study the Rees algebra of a link $I$ of a
Cohen--Macaulay prime ideal ${\frak p}$ still requires 
that $R$ be a Gorenstein ring. We carry out in Section $3$ the
determination of the canonical module of $R[It]$ in full detail. 
Self--linked prime ideals will correspond to the case in which the
associated graded ring is Gorenstein.

In the final section, we establish two multiplicity formulae. The
first says that the multiplicity of $R[Jt]$ and $R[It]$ are the same. 
Since there are explicit formulae for the multiplicities
of some Rees algebras of complete intersection (see \cite{HTU}), this
is useful for computational purposes. The second is an estimate
relating the multiplicity of a Cohen--Macaulay self-linked ideal
$I$ with Cohen--Macaulay first Koszul homology 
(e.g. ideals in the linkage class of complete intersections) to its
number of generators. 

\medskip

As a general reference for unexplained notations or basic results we
will use three invaluable sources: \cite{BH}, \cite{Mat} and
\cite{WVV}---the first one especially for the sections on the
canonical module and multiplicity. 

\medskip

Last but not least, it is a pleasure to thank our advisor, Professor
Wolmer V. Vasconcelos, for fruitful conversations and the referee for 
his valuable suggestions regarding this paper.

\section{Links of prime ideals in Cohen--Macaulay rings}
In this section we describe under which conditions the results of
\cite{CPV} can be extended to the link of a prime ideal in an
arbitrary Cohen--Macaulay ring. More precisely,  
Theorem~\ref{main} generalizes Theorem A to a
Cohen--Macaulay ring; the only assumption is that $R_{\frak p}$ cannot
be  a regular local ring. Actually, this is not a restriction because
the case of a regular ring was thoroughly studied in \cite{CPV}.

We first prove the theorem in the case of links of the maximal
ideal; we then show how the general case can be reduced to this
particular one. For that, one needs to know that if $I$ and $J$ are
linked ideals then ${\rm Ass}(R/I)\cup{\rm Ass}(R/J)={\rm Ass}(R/{\bf z})$;
in particular every associated prime of $I$ and $J$ has the same height.

An ideal-theoretic fact is needed in order to obtain the first result. 

\begin{Proposition}\label{propo1}
Let $A$ and $B$ be two ideals of a local ring $(R, {\frak m})$ and let
${\bf z}=z_1, ..., z_n$ be a regular sequence contained both in $A$
and in $B$. If $AB\subset ({\bf z})$ but $AB\not\subset {\frak
m}({\bf z})$ then $A$ and $B$ are both generated by regular
sequences of length $n$.
\end{Proposition}

\demo
Consider first the case $n=1$ and write $z=z_1$, for sake of simplicity.
By assumption, there exist $a\in A$ and $b\in B$ such that $ab=\alpha
z$ with $\alpha\not\in {\frak m}$. We may even assume $z=ab$. 
Pick now any $x$ in $A$; since $xb\in AB\subset (z)$, $xb=cz$ for some
$c\in R$ and so  
$$
xz=x(ab)=a(xb)=a(cz)=(ac)z.
$$
But $z$ is a regular element, thus $x=ac$, i.e. $x\in(a)$. Since $x$ is
an arbitrary element of $A$ we then conclude that $A=(a)$; similarly,
$B=(b)$.
Note that if the product of two elements is regular then they are both
regular, and this establishes the first case. 

Consider now the case $n>1$. Since $AB\not\subset {\frak m}({\bf z})$,
there exist $a\in A$ and $b\in B$ such that $ab=\alpha_1z_1+ \cdots +
\alpha_nz_n$ and at least one of the $\alpha_i$ does not belong to
${\frak m}$; by reordering the $z_i$'s we may assume $\alpha_n\not
\in{\frak m}$. Let ``${}^{-}$'' denote the homomorphic image modulo
$(z_1, ..., z_{n-1})$: we are then in the case of two ideals
$\overline{A}$ and $\overline{B}$ of a local ring $(\overline{R},
\overline{{\frak m}})$ both containing the regular element
$\overline{z}_n$. In addition, $\overline{A}\overline{B}\subset
(\overline{z}_n)$ but $\overline{A}\overline{B}\not\subset {\frak m}
(\overline{z}_n)$. Using the previous argument, we then conclude that
$\overline{A}=(\overline{a})$ and $\overline{B}=(\overline{b})$, 
which implies
$$
A=(z_1, ..., z_{n-1}, a)\quad {\rm and} \quad B=(z_1, ..., z_{n-1}, b),
$$
as claimed.\QED

\begin{Theorem}\label{Thm1}
Let $(R, {\frak m})$ be a $d$-dimensional Cohen--Macaulay local ring, 
which is not a regular local ring.  
Let $I=J\colon {\frak m}$ where $J=(z_1, ..., z_d)$ is a system of
parameters. Then $I^2=JI$.
\end{Theorem}
\demo
We divide the proof into two parts, first showing  that
${\frak m}I={\frak m}J$, and then that $I^2=JI$.

\medskip

Clearly, ${\frak m}I\supseteq {\frak m}J$; since $I=J\colon 
{\frak m}$  we have that ${\frak m}I\subset J$, so if ${\frak
m}I\not\subseteq {\frak m}J$, Proposition~\ref{propo1} with $A=I$,
$B={\frak m}$ and $({\bf z})=J$ would imply that ${\frak m}$ is
generated by a regular sequence, contradicting the assumption that $R$ 
is not a regular local ring. 

\medskip

It will suffice to show that the product of any two elements of $I$ is
contained in $JI$. Pick $a, b\in I$; since $ab\in I^2\subset {\frak
m}I\subset J$, $ab$ can be written in the following way 
\begin{equation}\label{eqT1}
ab=\sum_{i=1}^d \alpha_iz_i.
\end{equation}
The proof will be complete once it is shown that $\alpha_i\in
I$, for $i=1, ..., d$; to this end it is enough to prove that
$\alpha_ix\in J$ for any $x\in{\frak m}$. 

Since $R$ is a local ring regular sequences permute, thus
we can key on any of the coefficients, say $\alpha_d$. 
Note that $xb\in{\frak m}I={\frak m}J$, hence $xb$ satisfies 
an equation of the form
\begin{equation}\label{eqT2}
xb=\sum_{i=1}^d \beta_iz_i, \quad \quad \beta_i\in{\frak m}.
\end{equation}
Multiply equation (\ref{eqT1}) by $x$ and equation (\ref{eqT2})
by $a$; a quick comparison between these new equations yields the
following result
$$
\sum_{i=1}^d \alpha_i x z_i=\sum_{i=1}^d a\beta_i z_i,
$$
which can also be written as
\begin{equation}\label{eqT3}
\sum_{i=1}^d (\alpha_i x - a\beta_i) z_i=0.
\end{equation}
Equation (\ref{eqT3}) modulo the ideal $J_1=(z_1, ...,
z_{d-1})$ gives that $(\alpha_d x-a\beta_d)z_d$ is zero modulo $J_1$;
but $z_d$ is regular in $R/J_1$ hence
$$
\alpha_d x-a\beta_d\in J_1=(z_1, ..., z_{d-1})\subset J.
$$
Since $a\beta_d\in I{\frak m}\subset J$ we then have that
$\alpha_d x\in J$ as well, or equivalently $\alpha_d\in I$. \QED

\begin{Theorem}\label{main}
Let $(R,{\frak m})$ be a Cohen--Macaulay local ring and let ${\frak
p}$ be a prime ideal of $R$ of height $g$ such that $R_{\frak p}$ is not
a regular local ring. If $J=(z_1, ..., z_g)$ denotes a regular
sequence inside ${\frak p}$, then the link $I=J\colon {\frak p}$ is
equimultiple with reduction number $1$.
\end{Theorem}
\demo
It is enough to reduce to the case in which $(R, {\frak p})$ is a
local ring and then apply Theorem~\ref{Thm1}.
 
In order to show the equality $I^2=JI$ it suffices to
establish it at the associated primes of $R/JI$. From the short
exact sequence
$$
0 \rightarrow J/JI \longrightarrow R/JI \longrightarrow R/J
\rightarrow 0
$$
and the fact that $J/JI=J/J^2\otimes R/I=(R/I)^g$, it follows that 
the associated primes of $R/JI$ are contained in those of $R/J$ and of 
$R/I$. On the other hand, we claim that ${\rm Ass}(R/I)\subset {\rm
Ass}(R/J)$; indeed, $I=J\colon {\frak p}$ with ${\frak p}$ a prime
ideal implies that $I$ and ${\frak p}$ are directly linked.
Thus the associated primes of $R/JI$ must have height $g$, since $J$ is a
complete intersection. 

If we consider ${\frak q}\in{\rm Ass}(R/J)$, then either ${\frak q}={\frak
p}$ or ${\frak q}\not\supset {\frak p}$. The first case is the
situation of a link with the maximal ideal; in the second case
$I_{\frak q}=J_{\frak q}$ and so we are done. \QED 

\medskip

The next corollary generalizes Theorem B and Corollary C to the case
of Cohen--Macaulay links of prime ideals in arbitrary Cohen--Macaulay
rings. 

\begin{Corollary}
With the same assumptions of {\rm Theorem~\ref{main}} if $I$ is a
Cohen--Macaulay ideal then the associated graded ring of $I$ is
Cohen--Macaulay.
If in addition ${\rm dim}(R_{\frak p})\geq 2$ then the Rees algebra of
$I$ is Cohen--Macaulay as well.
\end{Corollary}
\demo
The proof is exactly the same as the one in \cite{CPV}. There, the
Cohen--Macaulayness of $I$ was a consequence of linkage in
Gorenstein rings; here, we have to ask explicitly for
it. 
\QED

\begin{Remark}{\rm 
The previous assertions leave out the cases in which $R_{\frak p}$ is
a regular local ring. When $R$ is Gorenstein this is taken care of by
\cite{CPV} provided that condition L${}_2$ is satisfied. We would like
to remark however that the associated graded ring of the link $I$ is
Cohen--Macaulay without any assumption if ${\frak p}$ is supposed to
be a Gorenstein ideal.
This is an easy consequence of Theorems A and B together with the fact
that in a Gorenstein setting a strongly Cohen--Macaulay ideal
satisfying property ${\cal F}_1$ has Cohen--Macaulay associated graded
ring (see \cite[Theorem 5.3]{HSV1}).   }
\end{Remark}

\begin{Remark}\label{rem26}{\rm
In the situation of Theorem~\ref{Thm1}
$$
I/J={\rm socle}(R/J) ={\rm Hom}_{R/J}(R/{\frak m},R/J)\simeq (R/{\frak
m})^s,
$$
where $s$ is the Cohen--Macaulay type of $R$. 
This and the fact that ${\frak m}I={\frak m}J$ imply that the
generators of $J$ are among the minimal generators of $I$ and that $s$
extra minimal generators are needed to get the whole $I$, i.e. 
$$
I=(J, a_1, ..., a_s).
$$
In a diagram this can be seen as follows (note that the numbers
represent the length of each quotient):
$$
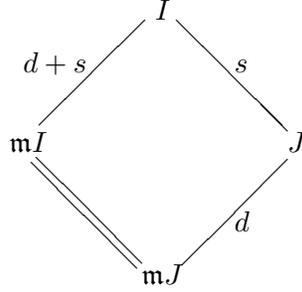
\begin{picture}(100,110)(-45,-55)
\put(0,50){$I$}
\put(50,0){$J$}
\put(-55,0){${\frak m}I$}
\put(-5, -50){${\frak m}J$}
\put(30,30){$s$}
\put(30,-30){$d$}
\put(-50,30){$d+s$}
\put(10, -43){\line(1, 1){40}}
\put(8, 50){\line(1, -1){42}}
\put(-44, 10){\line(1, 1){40}}
\put(-45, -2){\line(1, -1){40}}
\put(-47, -4){\line(1, -1){40}}
\end{picture}
$$ }
\end{Remark}

\begin{Remark}{\rm
Let $(R, {\frak m})$ be a Cohen--Macaulay local ring of dimension $d$ and
type $s\geq 2$ and let $I=J\colon {\frak m}$, with $J$ generated by a regular 
sequence of length $d$ inside ${\frak m}$. Then the following formula
is valid      
\begin{equation}\label{agg}
\lambda(\delta(I))={s+1 \choose 2},
\end{equation}
where $\delta(I)$ is the kernel of the natural surjection from
$S_2(I)$ to $I^2$ and $\lambda(\tratto)$ is the length function. 

\medskip

\demo
Let us consider the following commutative diagram
\begin{equation}\label{Diagramma1}
\begin{array}{ccccccccc}
& &0& &0& &0& & \\
& &\downarrow& &\downarrow& &\downarrow & &  \\
0 & \rightarrow & \star & \longrightarrow & \delta(I) &
\longrightarrow & \star\star & \rightarrow & 0 \\ 
& &\downarrow& &\downarrow& &\downarrow & &  \\
0 & \rightarrow & J\cdot I & \longrightarrow & S_2(I)
& \longrightarrow & S_2(I/J)& \rightarrow & 0 \\ 
& &\downarrow& &\downarrow& &\downarrow & &  \\
0 & \rightarrow & JI & \longrightarrow &I^2 &
\longrightarrow & I^2/JI& \rightarrow & 0 \\ 
& &\downarrow& &\downarrow& &\downarrow & &  \\
& &0& &0& &0& & 
\end{array}
\end{equation}
where $\cdot$ denotes the product in the symmetric algebra.
Since $J$ is a regular sequence inside $I$ we conclude that $J\cdot
I\simeq JI$ (see the proof of \cite[Theorem 2.4]{V1}) or,
equivalently, $\star=0$.   
Thus, using the additivity of the length on the third column of
(\ref{Diagramma1}) and the fact that $\delta(I)\simeq \star\star$,
one gets (see also \cite[Remark 2.4]{CPV})
$$
\lambda(I^2/JI)+\lambda(\delta(I))=\lambda(S_2(I/J)).
$$
Finally, (\ref{agg}) follows from $I/J\simeq (R/{\frak m})^s$
and $I^2=JI$. \QED }
\end{Remark}

\begin{Example}{\rm
The assumption of $R$ being Cohen--Macaulay cannot be dropped in
Theorem~\ref{main} as the following example shows; this
example, in particular, proves that the previous results cannot even be
extended to the case of systems of parameters (which are $d$-sequences) in
Buchsbaum rings.  

Let $k$ be any field, set $R=k[\![X, Y]\!]/(X^2, XY)=k[\![ x,
y]\!]$ where $x$ and $y$ denote the images of $X$ and $Y$ modulo
$(X^2, XY)$. It can be easily seen that $R$ is a one dimensional
Buchsbaum ring which is not Cohen--Macaulay (e.g. \cite[Example
5]{StVo}). Let $J=(y^3)$ be a system of parameters in $R$ and
consider $I=J\colon {\frak m}$, where ${\frak m}=(x, y)$; a
computation using the computer system {\em Macaulay} shows that
$I^2\not= JI$. }  
\end{Example}

\section{The canonical module} 
If $(R, {\frak m})$ is a Cohen--Macaulay local ring, a faithful
maximal Cohen--Macaulay module $\omega_R$ of type $1$ is 
called a {\em canonical module} of $R$. When it exists this module is
unique up to isomorphism. If $R$ is a Gorenstein local ring then
$\omega_R$ is isomorphic to $R$; moreover, if $\varphi\colon S
\longrightarrow R$ is a local homomorphism of Cohen--Macaulay local
rings such that $R$ is a finite $S$-module and $\omega_S$ is the
canonical module of $S$ then $\omega_R$ exists and is given by
$\omega_R\simeq {\rm Ext}^g_S(R, \omega_S)$, where $g={\rm
dim}(S)-{\rm dim}(R)$. 

In this section we describe the structure of the canonical
module of the Rees algebra of a Cohen--Macaulay, equimultiple ideal
$I$ with reduction number 1. This allows us to give a necessary and
sufficient condition for the associated graded ring of $I$ to be
Gorenstein.
As an application, we compute the canonical module of the
Rees algebra of the link $I$ of a prime ideal ${\frak p}$ satisfying
either condition L${}_1$ or L${}_2$. 

For convenience we recall some additional facts from linkage theory. 
Let $R$ be a Gorenstein local ring, $I$ an unmixed ideal of
grade $g$ and ${\bf z}$ a regular sequence of length $g$ inside $I$.
If we set $J={\bf z}\colon I$ then $I\sim J$. In addition, if 
$I\sim J$ then $I$ is Cohen--Macaulay if and only if $J$ is
Cohen--Macaulay. Furthermore, $\omega_{R/I}\simeq J/{\bf z}$ and
$\omega_{R/J}\simeq I/{\bf z}$. 

\begin{Theorem}\label{nuovo}
Let $R$ be a $d$-dimensional Gorenstein local ring and let $I$ be a
Cohen--Macaulay ideal of $R$, of height $g\geq 2$. Let $J=(a_1, ..., a_g) 
\subset I$ be a reduction of $I$ with $I^2=JI$. 
Then the canonical module of $R[It]$ has the form
$$
\omega_{R[It]}=(t(1, t)^{g-3}+Lt^{g-1}),
$$ 
where $L$ is given by $J\colon I$.
\end{Theorem}
\medskip
\demo
Set $A=R[Jt]$ and $B=R[It]$ and let $L=J\colon I$.
Since $I$ is a Cohen--Macaulay ideal, then $I=J \colon L$ and $L$ is a 
Cohen--Macaulay ideal (see \cite{PS}).

The canonical module of $A$ is the fractional ideal of $A$
given by $\omega_A=(1,t)^{g-2}[-1]$ (see \cite{B, HV}). 
More explicitly 
$$
\omega_A=Rt+Rt^2+\cdots +Rt^{g-1}+Jt^g+J^2t^{g+1}+\cdots .
$$
Since $J$ is a reduction of $I$, the extension 
$A \hookrightarrow B$ is finite; therefore the canonical module of $B$ 
is given by  
\begin{equation}\label{localize}
\omega_B = {\rm Hom}_A(B,\omega_A);
\end{equation}
as a graded ideal, which we represent as
$$
\omega_B = \omega_1t + \omega_2t^2+ \cdots .
$$
Indeed, localizing $\omega_B$ at the multiplicative set given by 
the powers of $x$, where $x\in J$ is a regular element, we get 
that $\omega_{B_x}=\omega_{R_x}[t][-1]=R_x[t][-1]$, since $B_x=R_x[t]$
and $R_x$ is Gorenstein. 
This means that $1$ is the initial degree of $\omega_{B_x}$.
However, as $x$ is regular on $B$ (and therefore on $\omega_B$ as
well) we have the inclusion $\omega_0 \hookrightarrow
(\omega_0)_x =0$, so that $\omega_0=0$.  

The problem is to recover the $\omega_i$ from the formal expression
of $\omega_B$ as the dual above. 
Consider the natural exact sequence of $R[Jt]$-modules,
\begin{equation}
0 \rightarrow  A \longrightarrow B \longrightarrow C \rightarrow 0,
\label{eq:1}
\end{equation}
where $C$ is $I/J+I^2/J^2+\cdots+I^n/J^n+\cdots$. Note that $C$ 
is a Cohen--Macaulay $R[Jt]$-module of dimension $d$.
Indeed, from (\ref{eq:1}) and the depth lemma (see 
\cite[Lemma 1.2.9]{BH}) we have that 
${\rm depth}(C)\geq d$. On the other hand, the dimension of $C$ 
is less or equal than $d$ since $JR[Jt]$ has height $1$ and is
contained in the annihilator of $C$. 

Dualizing (\ref{eq:1}) with $\omega_A$, we have the following exact
sequence 
\begin{equation}
0\rightarrow \omega_B \longrightarrow \omega_A \longrightarrow
\omega_C={\rm Ext}^1_A(C,\omega_A) \rightarrow 0. 
\label{eq:0}
\end{equation}
It follows that if one knows what ${\rm Ext}^1_A(C,\omega_A)$ looks like, 
$\omega_B$ can be recovered.
In order to describe ${\rm Ext}^1_A(C,\omega_A)$ we have to examine
more closely $C$.
Since $I^2=JI$, we have that $C=I/J\cdot R[Jt][-1]$,
i.e. $C_i = (I/J\cdot R[Jt])_{i-1}$;
therefore from the presentation $R[x_1,...,x_g] \longrightarrow R[Jt]
\rightarrow 0$, we have the following surjection
\begin{equation}
I/J\otimes R[x_1,...,x_g][-1]=I/J[x_1,...,x_g][-1] \longrightarrow
C=I/J\cdot R[Jt][-1] \rightarrow 0. 
\label{eq:2}
\end{equation}
We want to show that 
$$
C\simeq I/J[x_1,...,x_g][-1].
$$
Let $K$ be the kernel of the map in (\ref{eq:2})
\begin{equation}
0 \rightarrow  K \longrightarrow I/J[x_1,...,x_g][-1] \longrightarrow
 C \rightarrow 0;
\label{eq:3}
\end{equation}
to show that $K=0$ it will be enough to
show that $K_{\frak q}=0$ for all ${\frak q}\in {\rm Ass}(I/J[x_1,
..., x_g])$.
Note that $I/J=(J\colon L)/J=\omega_{R/L}$ is a Cohen--Macaulay module
of dimension $d-g$. In particular $I/J[x_1,..., x_g]$ is a
Cohen-Macaulay module of dimension $d$. By the depth lemma we also
have that ${\rm depth}(K)\geq d$, hence $K$ is Cohen-Macaulay as well.
Next, we observe that the associated primes of
$K$ are included in the associated primes of $I/J[x_1,..,x_g]$;
and that ${\rm Ass}(I/J)={\rm Ass}(R/L)$.
Localize (\ref{eq:3}) at ${\frak q}$ to get the exact sequence
\begin{equation}\label{nove}
0 \rightarrow K_{\frak q} \longrightarrow (I/J)_{\frak q}[x_1,
...,x_g][-1]=\omega_{R_{\frak q}/L_{\frak q}}[x_1, ...,x_g][-1]
\longrightarrow C_{\frak q} \rightarrow 0, 
\end{equation}
of Cohen--Macaulay modules of dimension $g$. 
Let $D=R_{\frak q}/L_{\frak q}$; $D$ is an Artinian ring with maximal
ideal ${\frak M}={\frak q}R_{\frak q}/L_{\frak q}$.
Set $D[x_1, ...,x_g]=D[\underline{x}]$ and dualize (\ref{nove}) with 
${\rm Hom}_{D[\underline{x}]}(\tratto,\omega_{D[\underline{x}]})$ to 
obtain the exact sequence 
$$
0 \rightarrow {\rm Hom}_{D[\underline{x}]} (C_{\frak q},
\omega_{D[\underline{x}]}) \longrightarrow
{\rm Hom}_{D[\underline{x}]}(\omega_{D}[\underline{x}][-1],
\omega_{D[\underline{x}]}) \longrightarrow 
{\rm Hom}_{D[\underline{x}]}(K_{\frak q}, \omega_{D[\underline{x}]}) 
\rightarrow 0 
$$
or
\begin{equation}
0 \rightarrow C_{\frak q}^{\vee} \longrightarrow
D[\underline{x}][1-g] \longrightarrow K_{\frak q}^{\vee} \rightarrow 0,
\end{equation}
of Cohen--Macaulay modules of dimension $g$. It will be enough
to prove that ${\rm height}(C_{\frak q}^{\vee})\geq 1$, since this
would imply that ${\rm dim}(K_{\frak q}^{\vee})\leq g-1$. 
In order to show ${\rm height}(C_{\frak q}^{\vee})\not=0$ we
have to show that $C_{\frak q}^{\vee}\not\subseteq {\frak
M}[\underline{x}]$, as ${\frak M}[\underline{x}]$ is the unique 
minimal prime of $D[\underline{x}]$. Assume otherwise; thus
${\rm Ann}_D(C_{\frak q}^{\vee})\supseteq
{\rm Ann}_D({\frak M}[\underline{x}])\supseteq {\frak
M}^{r-1}\not= 0$, where $r$ is the index of nilpotency of
the maximal ideal ${\frak M}$ of the Artinian ring $D$.
This leads to a contradiction since
${\rm Ann}_D(C_{\frak q}^{\vee})=0$. Indeed,
$C_{\frak q}^{\vee}={\rm
Hom}_{D[\underline{x}]}(C_{\frak q}, \omega_{D[\underline{x}]})$ and 
$C_{\frak q}\simeq {\rm Hom}_{D[\underline{x}]}(C_{\frak
q}^{\vee}, \omega_{D[\underline{x}]})$ imply that ${\rm Ann}_{D}
(C_{\frak q}^{\vee})={\rm Ann}_{D}(C_{\frak q})\subseteq 
{\rm Ann}_{D}((I/J)_{\frak q})=0$.
The last part follows from the fact that $C_{\frak q}= (I/J)_{\frak q}
+ (I^2/J^2)_{\frak q}+ \cdots$ and that $(I/J)_{\frak
q}=\omega_{D}$ is $D$-faithful, since $D={\rm End}_{D}(\omega_D)$.
Hence
$$
C=\omega_{R/L} [x_1,...,x_g][-1]=\omega_{R/L [x_1,...,x_g]}[g-1].
$$
Finally, by duality we have that
$$
\omega_C={\rm Ext}^1_A(C,\omega_A)={\rm
Ext}^1_A(\omega_{R/L[x_1,...,x_g]}[g-1],\omega_A)=
R/L[x_1,...,x_g][-g+1].  
$$
Notice that $R/L[x_1, ..., x_g]\simeq {\rm gr}_J(R)\otimes
R/L$; therefore two of the modules in the short exact sequence (\ref{eq:0})
are
$$
\begin{array}{l}
\omega_A=(t(1, t)^{g-2}), \\
\omega_C=({\rm gr}_J(R)\otimes R/L)[-g+1].
\end{array}
$$
From (\ref{eq:0}) we conclude that $\omega_i\simeq R$ for all $i=1, ...,
g-2$ and that $\omega_{g-1}$ is given by $R/\omega_{g-1}\simeq R/L$.
If we show that $\omega_{g-1}\subseteq L$ we get that  
$\omega_{g-1}\simeq L$. In order to show such an inclusion,
we have to look at the structure of the canonical module in a
different way.
This alternative approach is based on the following observation
$$
\omega_B={\rm Hom}_A(B,\omega_A)\simeq \omega_A\colon_{{}_{R[t]}} B.
$$
Thus $\omega_{g-1}$ is going to satisfy the following set of
relations
$$
\omega_{g-1} R\subset R, \ \omega_{g-1} I\subset J, \ \omega_{g-1}
I^2\subset J^2, \ldots .
$$
The first inclusion says that $\omega_{g-1}\subseteq R$; since
$I^2=JI$, the only important condition is given by
$\omega_{g-1} I \subset J$. From this equation we get the desired
inclusion since $I=J\colon L$. 

Finally, observe that $\omega_A$ modulo the submodule of $\omega_B$ 
generated by the elements up to degree $g-1$ is already $\omega_C$. \QED

\medskip

\begin{Corollary}\label{nuovo2}
Let $I$ be an ideal as in {\rm Theorem~\ref{nuovo}}, then ${\rm
gr}_I(R)$ is Gorenstein if and only if $I=J\colon I$.
\end{Corollary}
\demo
By a result of Herzog--Simis--Vasconcelos (see
\cite[Corollary 2.5]{HSV}), ${\rm gr}_I(R)$ is Gorenstein if and only if the
canonical module of $R[It]$ has the so called {\em expected form}, i.e.
$\omega_{R[It]}\simeq \omega_R t (1, t)^m$ for some $m\geq -1$.
By Theorem~\ref{nuovo}, 
$$
\omega_{R[It]}=(t(1, t)^{g-3}+Lt^{g-1}),
$$ 
hence it assumes the expected form if and only if $J\colon I=L=I$. 
\QED

\medskip

Let us now apply Theorem~\ref{nuovo} and Corollary~\ref{nuovo2} to the
{\em leitmotiv} of this paper, i.e. to the case of links of prime ideals.

\begin{Corollary}\label{Thmcano}
Let $R$ be a Gorenstein local ring of dimension $d$, ${\frak p}$
a Cohen--Macaulay prime ideal of codimension $g\geq 2$, $J=(x_1, ...,
x_g)$ a complete intersection contained in ${\frak p}$, and let $I=J\colon
{\frak p}$. If either {\rm L}${}_1$ or {\rm L}${}_2$ holds, then the
canonical module of $R[It]$ has the form 
$$
\omega_{R[It]}=(t(1, t)^{g-3}+{\frak p}t^{g-1}).
$$
\end{Corollary}
\medskip
\demo
From \cite{CPV}, $I$ is an equimultiple ideal with reduction number $1$ and
Cohen--Macaulay Rees algebra, so that Theorem~\ref{nuovo} applies.
\QED

\medskip

\begin{Remark}{\rm 
In the case of Corollary~\ref{Thmcano} the isomorphism 
$$
C\simeq I/J[x_1, ..., x_g][-1]
$$
can be established in a more direct manner. Indeed, the only associated 
prime of $I/J$ is ${\frak p}$. Localizing (\ref{eq:3}) at ${\frak p}$, 
we get the short exact sequence 
$$
0 \rightarrow K_{{\frak p}} \longrightarrow 
(I/J)_{{\frak p}}[x_1,...,x_g][-1] \longrightarrow C_{{\frak p}}
\rightarrow 0, 
$$
where $(I/J)_{{\frak p}}[x_1,..,x_g]$ is isomorphic to a polynomial 
ring of dimension $g$ over a field and $K_{{\frak p}}$ to one of
its ideals. Thus ${\rm dim}(C_{{\frak p}})\leq g-1$. 
But this is a contradiction since $C_{{\frak p}}$ is still a 
Cohen--Macaulay module of dimension $g$.} 
\end{Remark}

\section{On the multiplicity of special rings}
In this section we study the multiplicity of
two families of rings: Rees algebras of links of the maximal ideal of
a Gorenstein local ring, and quotients of a Gorenstein local ring
modulo a self-linked ideal.

\subsection{The Rees algebra of links} 
The following result is a very useful tool in the computation of the 
multiplicity of certain Rees algebras.

\begin{Theorem}\label{multe}
Let $(R,{\frak m})$ be a $ d$-dimensional Gorenstein local ring and
let $I=J\colon{\frak m}$ where $J$ is a regular sequence of $R$ of
length $d$. Let ${\frak M}=({\frak m}, It)$ be the maximal  
homogeneous ideal of $R[It]$ and ${\frak N}=({\frak m}, Jt)$ be the maximal 
homogeneous ideal of $R[Jt]$. If either {\rm L}${}_1$ or {\rm L}${}_2$
holds, then the multiplicities of those two Rees algebras with respect
to their maximal ideals are equal,
$$
e({\frak M}, R[It])=e({\frak N}, R[Jt]).
$$
\end{Theorem}
\demo
$R[It]$ is a finite $R[Jt]$-module of rank $1$, as $J$ is a 
reduction of $I$. Applying \cite[Corollary 4.6.9]{BH} 
to the ring $(R[Jt], {\frak N})$ and to the module $M=R[It]$ we get
that 
$$
e({\frak N}R[It], R[It])=e({\frak N}, R[Jt]).
$$
On the other hand, since $I=J \colon {\frak m}$, or more specifically
${\frak m}I\subset J$, and $I^2=JI$ (with $J$ a minimal reduction of $I$)
we have that ${\frak m}I={\frak m}J$. So, the following calculation
goes through
\begin{eqnarray*}
({\frak N}R[It])\, {\frak M} 
                      &=& ({\frak m}+Jt+JIt^2+JI^2t^3+\cdots) ({\frak
                          m}+It+I^2t^2+I^3t^3+\cdots)  \\ 
                      &=& {\frak m}^2+{\frak m}It+JIt^2+JI^2t^3+\cdots  \\
                      &=& {\frak M}^2,
\end{eqnarray*}
that is to say ${\frak N}R[It]$ is a reduction of ${\frak M}$. 
Therefore the assertion follows from \cite[Lemma 4.5.5]{BH} when applied 
to the ring $(R[It], {\frak M})$ and to the module $M=R[It]$. 
\QED

\medskip

\begin{Example}{\rm 
Thanks to Theorem~\ref{multe}, the problem of computing $e({\frak M},
R[It])$ is then reduced to the case in which we are dealing with the
Rees algebra $R[Jt]$ of a regular sequence. The latter problem has
been recently studied in \cite{HTU} in the case in which $R$ is an
homogeneous algebra over a field $k$ and $J$ is an ${\frak m}$-primary
ideal generated by an homogeneous regular sequence. 

The next example is taken from \cite{V}.
Let $(R, {\frak m})$ be a $3$-dimensional regular local ring
containing a field $k$ and assume that the maximal ideal is given by
the following regular system of parameters ${\frak m}=(x, y, z)$. Let
$I$ and $J$ be the following ideals   
$$
I=J\colon{\frak m}\quad \quad J=(x, y^2, z^2).
$$
It turns out that $I=(x, y^2, yz, z^2)$; moreover, since the generators
of $J$ are homogeneous of degrees
$$
a_1={\rm deg}(x)=1, \quad a_2={\rm deg}(y^2)=2, \quad a_3={\rm deg}(z^2)=2,
$$
our result combined with \cite[Corollary 1.5]{HTU} gives directly
$$
e({\frak M}, R[It])=e({\frak N}, R[Jt])=(1+a_1+a_1a_2)e(R)=(1+1+1\cdot
2)e(R)=4\cdot 1=4. 
$$}
\end{Example}

\subsection{Quotients of self-linked primes}
Finally, we study the multiplicity of $ R/I$ where $I$ is a
self-linked ideal. 
We give a lower  bound for this number when, essentially, the
first Koszul homology module $H_1(I)$ is Cohen--Macaulay.

\begin{Theorem}
Let $R$ be a Gorenstein local ring and $I$ a self-linked 
ideal of $R$ of height $g$; assume that $I$ is generically
Gorenstein and that it lies in the linkage class of a complete
intersection {\rm (}licci{\rm )}. 
Let $\beta_1$ denote the minimal number of generators of $I$.
Then 
$$
e(R/I)\geq {\beta_1-g+1 \choose 2},
$$
where $e(R/I)$ denotes the multiplicity of $R/I$.
\end{Theorem}
\demo
Write $I=J\colon I$, where $J$ is an ideal generated by a regular
sequence of length $g$ inside $I$. Obviously $I/J$ is an $R/I$-module.
For every associated prime ${\frak q}$ of $I$, $(I/J)_{\frak q}\simeq
\omega_{(R/I)_{\frak q}}\simeq (R/I)_{\frak q}$, which shows that
$I/J$ is an $R/I$-module of rank $1$, and hence so is $S_2(I/J)$. But
the latter module is Cohen--Macaulay since 
$$
{\rm depth}(H_1(I))={\rm depth}(S_2(I/J))={\rm
depth}(S_2(\omega_{R/I}))  
$$
(see \cite[Theorem 3.1]{Ul}) and since $H_1(I)$ is Cohen--Macaulay by
\cite{Hun}. Finally, from a general formula due to D. Rees
(see \cite[Corollary 5.3.3]{Sa}), we have that
$$
e(R/I)\geq \mu(S_2(I/J))\geq {\beta_1-g+1\choose 2},
$$
as claimed.
\QED

\end{document}